\documentclass[12pt, reqno]{amsart}
\usepackage{amsfonts}
\usepackage{ifthen}
\usepackage{amsthm}
\usepackage{amsmath}
\usepackage{graphicx}
\usepackage{amscd,amssymb,amsthm}
\newtheorem{problem}{Problem}[section]\newtheorem{proposition}{Proposition}[section]

\newcounter{minutes}
\divide\time by 60
\newcounter{hours}
\multiply\time by 60 \addtocounter{minutes}{-\time}

\newcommand{\real}{\operatorname{Re}}

\newtheorem{theorem}{Theorem}[section]
\newtheorem{lemma}{Lemma}[section]

\title{On the geometric properties of the Bessel-Struve kernel function}

\author[Saiful R. Mondal]{Saiful R. Mondal$^\ast$}
\address{Department of Mathematics and Statistics, Collage of Science,
King Faisal University, Al-Hasa 31982, Hofuf, Saudi Arabia.}
\email{smondal@kfu.edu.sa}
\thanks{$\ast$ The author thanks the  Deanship of Scientific Research at King Faisal University for funding his work under project number 150244.}


\begin{document}

\begin{abstract}
This paper introduce the  Bessel-Struve kernel functions $\mathcal{B}_\nu$ defined on the unit disk in the complex plane. We studies the close-to-convexity of $\mathcal{B}_\nu$ with respect to several starlike functions. Sufficient condition on $\nu$ for which the function $z\mathcal{B}_\nu$ is starlike is given.
\end{abstract}
\maketitle
\section{Introduction and Preliminaries}
\subsection{Bessel-Struve Kernel functions}

The function $J_\nu$, known as the Bessel function of the first kind of order $\nu$, is  a particular solution of the second order Bessel differential equation 
\begin{align*}
z^2y''(z)+zy'(z)+(z^2-\nu^2)y(z)=0.
\end{align*}
This function has the power series representation
\begin{align}\label{eqn:Bessel}
J_\nu(z)= \sum_{n=0}^\infty \frac{(-1)^n}{n! \Gamma{(\nu+n+1)} }\left(\frac{z}{2}\right)^{2n+\nu}, \quad |z|<\infty.
\end{align}
On the other hand the modified Bessel function $I_\nu(z)$ is the particular solution of the differential equation
\begin{align*}
z^2y''(z)+zy'(z)-(z^2-\nu^2)y(z)=0,
\end{align*}
and have the series representation
 \begin{align}\label{eqn:modified-Bessel}
I_\nu(z)= \sum_{n=0}^\infty \frac{1}{n! \Gamma{(\nu+n+1) }}\left(\frac{z}{2}\right)^{2n+\nu}, \quad |z|<\infty.
\end{align}

The Struve function of order $\nu$ given by
\begin{equation}\label{Strv}
\mathtt{H}_\nu(z):=  \sum_{n=0}^\infty \frac{(-1)^k}{\mathrm{\Gamma}{\left(n+\nu+\frac{3}{2}\right)} \mathrm{\Gamma}{\left(n+\frac{3}{2}\right)}} \left(\frac{z}{2}\right)^{2n+\nu+1}
\end{equation}
is a particular solution of the non-homogeneous Bessel differential equation
 \begin{align}\label{StrvDE}
 z^2 y''(z)+ z y'(z)+(z^2-\nu^2)y(z)= \frac{4\left(\frac{z}{2}\right)^{\nu+1}}{\sqrt{\pi}\; \Gamma{\left(\nu+\frac{1}{2}\right)}}.
 \end{align}
 A solution of the non-homogeneous modified Bessel equation
 \begin{align}\label{modStrvDE}
 z^2 y''(z)+ z y'(z)-(z^2+\nu^2)y(z)= \frac{4\left(\frac{z}{2}\right)^{\nu+1}}{\sqrt{\pi}\; \Gamma{\left(\nu+\frac{1}{2}\right)}}.
 \end{align}
 yields the modified Struve function
 \begin{align}\label{modStrv}\notag
\mathtt{L}_\nu(z)&:=- i e^{-\frac{i \nu \pi}{2}}\; \mathtt{H}_\nu(i z)\\&=  \sum_{n=0}^\infty \frac{1}{\mathrm{\Gamma}{\left( n +\nu+\frac{3}{2}\right)} \mathrm{\Gamma}{\left( n+\frac{3}{2}\right)}} \left(\frac{z}{2}\right)^{2n+\nu+1}.
\end{align}
The Struve functions occur in areas of physics and applied mathematics, for example, in water-wave and surface-wave problems \cite{Ahmadi-Widnall}, as well as in problems on unsteady aerodynamics \cite{Shaw}. The Struve functions are also important in particle quantum dynamical studies of spin decoherence \cite{shao} and nanotubes \cite{Pedersen}.

Consider the Bessel-Struve kernel function $\mathcal{B}_\nu$ defined on  the unit disk $\mathbb{D}=\{z: |z|<1\}$ as
\begin{align}
\mathcal{B}_\nu(z):= j_\nu(i z)- i h_\nu(iz), \quad \nu> -\frac{1}{2},
\end{align}
where, $j_\nu(z):= 2^\nu z^{-\nu} \Gamma{(\nu+1)} \mathtt{J}_\nu(z)$ and $h_\nu(z):= 2^\nu z^{-\nu} \Gamma{(\nu+1)} \mathtt{H}_\nu(z)$ are respectively known as the normalized Bessel functions and the normalized Struve functions of first kind of index $\nu$. The Bessel-Struve transformation and Bessel-Struve kernel functions  are appeared in many article. In \cite{Hamem}, Hamem et. al. studies an analogue of the Cowling's Price theorem for the Bessel–Struve transform defined on real domain and also provide Hardy type theorem associated with this transform. The Bessel-Struve intertwining operator on $\mathbb{C}$ is considered in \cite{GS, KS}. The fock space of the Bessel-Struve kernel functions is discussed in \cite{GS2}.
The kernel $z \mapsto \mathcal{B}_\nu(\gamma z)$, $\gamma \in \mathbb{C}$ is the unique solution of the initial value problem
\begin{align}\label{eqn:Bessel-struve-opt}
\mathcal{L}_\nu u(z)= \lambda^2 u(z), \quad u(0)=1, u'(0)=\frac{\lambda \Gamma{(\nu+1)}}{\sqrt{\pi}\Gamma{(\nu+\frac{3}{2}})}.
\end{align}
Here $\mathcal{L}_\nu$ , $\nu >-1/2$ is the Bessel-Struve operator given by
\begin{align}\label{eqn:Bessel-struve-opt-1}
\mathcal{L}_\nu(u(z)):= \frac{d^2u}{dz^2}(z)+\frac{2\nu+1}{z}\left(\frac{du}{dz}(z)-\frac{du}{dz}(0)\right).
\end{align}
Now from \eqref{eqn:Bessel} and \eqref{modStrv}, it is evident that  $\mathcal{B}_\nu$ (taking $\gamma=1$)  possesses the power series
\begin{align}\label{eqn:B-S-power}
\mathcal{B}_{\nu}(z):= \sum_{n=0}^\infty \frac{ {\Gamma(\nu+1)}\Gamma{\left(\frac{n+1}{2}\right)}}{\sqrt{\pi} n! \Gamma\left(\frac{n}{2}+\nu+1\right)} z^n.
\end{align}
The kernel $\mathcal{B}_\nu$ also have the integral representation
\begin{align}\label{kernel-intgra-rep}
\mathcal{B}_\nu( z):=\frac{ 2 \Gamma{(\nu+1)}}{\sqrt{\pi}\Gamma{\left(\nu+\frac{1}{2}\right)}}\int_{0}^1 (1-t^2)^{\nu-\frac{1}{2}} e^{ zt} dt.
\end{align}
The  identity \eqref{eqn:Bessel-struve-opt} and \eqref{eqn:Bessel-struve-opt-1} together imply that $\mathcal{B}_\nu$ satisfy the differential equation
\begin{align}\label{eqn:kumar-hypr-ode}
z^2\mathtt{g}''_\nu(z)+(2\nu+1)z\mathtt{g}'_\nu(z)-z\mathtt{g}_\nu(z)= z M ,
\end{align}
where $M=2\Gamma(\nu+1)\left(\sqrt{\pi}\;\Gamma(\nu +\frac{1}{2})\right)^{-1}$.

Another significance is that $\mathcal{B}_\nu$ can be express as the sum of the modified Bessel and the modified Struve function of first kind of order $\nu$. For the sake  of completeness, in the following result we established this relation.

\begin{proposition}\label{prop1}
For $\nu>0$, the following identity holds:
\begin{align*}
z^\nu \mathcal{B}_\nu(z)= 2^\nu \Gamma{(\nu+1)}\left(\mathtt{I}_\nu(z)+\mathtt{L}_\nu(z)\right).
\end{align*}
\end{proposition}
The function $\mathcal{B}_\nu$ have the following recurrence relation which is useful in sequel.
\begin{proposition}\label{prop2}
For $\nu > 0$, the following recurrence relation holds for $\mathcal{B}_\nu$:
\begin{align}\label{eqn-recurrence}
z \mathcal{B}'_\nu(z)= 2 \nu \mathcal{B}_{\nu-1}(z)-2 \nu \mathcal{B}_\nu(z).
\end{align}
\end{proposition}

\subsection{Starlike and close-to-convex functions}
Let $\mathbb{D}=\{z : |z|<1\}$ be the unit disk and $\mathcal{A}$ be the class of all analytic functions $f$ defined on $\mathbb{D}$ such that
$f(0)=0=f'(0)-1$. Clearly each $f \in \mathcal{A}$ have the power series
\begin{align}
f(z)= z+\sum_{n=2}^\infty a_n z^n.
\end{align}
 A function $f \in \mathcal{A}$ is said to be starlike if $f(\mathbb{D})$ is starlike with respect to the origin.  Now if for any starlike function $g$ and for some real number $\beta$, we have  $\real (e^{i \beta} f'(z)/g(z))>0$, then the function $f$ is said to be close-to-convex with respect to the starlike function $g$.  A function $f \in \mathcal{A}$ is convex if $f(\mathbb{D})$ is convex. The starlike and convex functions can be represent analytically  as
\[ \real\left(\frac{zf'(z)}{f(z)}\right)>0 \quad \text{and} \quad  \real\left(1+\frac{zf''(z)}{f'(z)}\right)>0,\]
 respectively. Traditionally the class of starlike functions is denoted as  $\mathcal{S}^{\ast}$, while the class of close-to-convex, and convex functions are denoted respectively as $\mathcal{C}$ and $\mathcal{K}$. These classes also be generalized by order $\lambda \in [0,1)$ with the analytical formulation as follows:
 \begin{align*}
   f\in \mathcal{S}^{\ast}(\lambda) &\Leftrightarrow \real\left(\frac{zf'(z)}{f(z)}\right)>\lambda, \\
   f\in \mathcal{C}(\lambda) &\Leftrightarrow \real\left(1+\frac{zf''(z)}{f'(z)}\right)>\lambda, \\
   f\in \mathcal{K}(\lambda) &\Leftrightarrow \real\left(\frac{zf'(z)}{g(z)}\right)>\lambda, \quad \text{for some}\quad  g\in \mathcal{S}^{\ast}.
 \end{align*}
 According to the Alexander duality theorem \cite{Alexander}, the function  $f : \mathbb{D} \to \mathbb{C}$ is in $\mathcal{C}(\nu)$ if and only if $z \to zf'(z)$ is starlike of order $\nu$. Here we remark that the definition of $\mathcal{C}(\nu)$ is also valid for non-normalized analytic function $f : \mathbb{D} \to \mathbb{C}$ with the property $f'(0)\neq 0$.

 Let  introduce another subclass of $\mathcal{S}^{\ast}(\lambda)$ consisting of functions $f$ for which
 \begin{align}\label{eqn:sc-s1-lambda}
  \left|\frac{zf'(z)}{f(z)}-1\right|< 1-\lambda,
 \end{align}
 and denoted the class as $\mathcal{S}_1(\lambda)$. The Alexander duality theorem can be apply to the class $\mathcal{S}_1(\lambda)$ and a function $f$ is said to be in the class $\mathcal{C}_1(\lambda)$ if $z f'(z) \in \mathcal{S}_1(\lambda)$.

 Following result is required in sequel.
 \begin{lemma}\cite{Owa-Srivastava}\label{lem:suf-starlike}
 Let $\lambda\in [0,1/2]$ be fixed and $\beta \geq 0$. If $f \in \mathcal{A}$ and
 \begin{align}\label{eqn:suf-starlike}
 \left|\tfrac{zf'(z)}{f(z)}-1\right|^{1-\beta}\left|\tfrac{zf''(z)}{f'(z)}\right|^\beta < (1-\lambda)^{1-2\beta}\left(1-\tfrac{3\lambda}{2}+\lambda^2\right)^\beta,
 \end{align}
 for all $z \in \mathbb{D}$, then $f \in \mathcal{S}^{\ast}(\lambda)$.
  \end{lemma}

Next we state our main results which are proved in Section $\ref{sec-2}$ by using Lemma \ref{lem:suf-starlike}.

 \begin{theorem}\label{thm:starlike}
Let the function $\mathcal{B}_\nu$ as defined in \eqref{eqn:Bessel-struve-opt} satisfy the inequality
\begin{align}\label{eqn:thm-starlike}
\left|\frac{z\mathcal{B}'_\nu(z)}{\mathcal{B}_\nu(z)}\right| <1-\lambda,
\end{align}
for $\lambda \in [0, 1/2]$. Then $z \mathcal{B}_\nu \in \mathcal{S}^\ast(\lambda)$.
\end{theorem}

Now we will introduce a subclass of $\mathcal{S}^{\ast}(\lambda)$ consisting of functions $f$ satisfying the inequality
\begin{align*}
\left|\frac{zf'(z)}{f(z)}-1\right|< 1-\lambda,
\end{align*}
is known as the class of the starlike functions with respect to 1 and denoted as $\mathcal{S}_1^{\ast}(\nu)$.

In our next result, we obtain sufficient condition by which the Bessel-Struve kernel functions is starlike with respect to $1$.

\begin{theorem}\label{thm:starlike1}
Let the function $\mathcal{B}_\nu$ as defined in \eqref{eqn:Bessel-struve-opt} satisfy the inequality
\begin{align}\label{eqn:thm-starlike}
\left|\frac{z\mathcal{B}''_\nu(z)}{\mathcal{B}'_\nu(z)}\right| <\frac{2-3\lambda+2\lambda^2}{2(1-\lambda)},
\end{align}
for $\lambda \in [0, 1/2]$. Then $\mathcal{B}_\nu \in \mathcal{S}^\ast_1(\lambda)$.
\end{theorem}

Following problem is  well known in the literature:
 \begin{problem}\label{prob-1}
 Find the conditions on $a_n$ such that
 \begin{align}\label{eqn:f}
 f(z)= z+\sum_{n=2}^\infty a_n z^n
 \end{align}
  is close-to-convex, or starlike or convex or any other subclasses of univalent functions.
 \end{problem}
Now in accordance with the Problem \ref{prob-1}, we need to find the condition on $\nu$ such that the Bessel-Struve kernel functions or it's normalization, will be in  any of the classes mention above.

 There are many results in the literature (see. \cite{Acharya, Ali-Lee-Mondal, MS, AS} and reference their in) which answer the above problem.  As per requirement for this work, we listed few of them here. Here we would like to remark that the functions
 \begin{align*}
 z, \frac{z}{1-z}, \frac{z}{1-z^2}, \frac{z}{(1-z)^2} \quad \text{and} \quad \frac{z}{1-z+z^2},
 \end{align*}
 and their particular rotations
 \begin{align*}
 \frac{z}{1+z}, \frac{z}{1+z^2}, \frac{z}{(1+z)^2} \quad \text{and} \quad \frac{z}{1+z+z^2}
 \end{align*}
 are the only nine functions which have integer coefficients and are starlike univalent in $\mathbb{D}$ (See \cite{BF}). The sufficient coefficient conditions for which a function  $f \in \mathcal{A}$ is close-to-convex can be easily obtain atleast when the corresponding starlike functions is one of the above listed form. In this article we will consider for $z$, $z/(1-z)$ and $z/(1-z^2)$.

\begin{lemma}\label{lemma-ac} \cite{Acharya}
Let $\{a_n\}_{n=1}^\infty$ be a sequence of non-negative real numbers such that $a_1=1$, $\underline{\Delta} a_n \geq 0$ when $n \geq 1$ and  $\underline{\Delta}^2 a_n$ when $n \geq 2$. Then the function $f$, defined in \eqref{eqn:f},  is starlike and close-to-convex with respect to the starlike functions $z$ and $z/(1-z)$.  Here $\underline{\Delta} a_n= n a_n- (n+1) a_{n+1}$ and
$\underline{\Delta}^{m+1} a_n=\underline{\Delta}^{m}(\underline{\Delta} a_n) $, $m=1, 2, \cdots$.
\end{lemma}
The starlikeness and close-to-convexity of  $z\mathcal{B}_\nu$  is obtained by using Lemma \ref{lemma-ac}  in the following result.
\begin{theorem}\label{theorem-1}
For $\nu \geq 1/2$,  the normalized Bessel-Struve kernel function $z\mathcal{B}_\nu$ is starlike. The function $z\mathcal{B}_\nu$ is also close-to-convex with respect to the starlike functions $z$ and $z/(1-z)$.
\end{theorem}

It can be observed that Theorem $\ref{theorem-1}$ can also be proved by using the following lemma given in \cite{MS}.

\begin{lemma}\cite{MS}
Let $\{a_n\}_{n \geq 1}$ be a sequence of positive real number such that $a_1 \geq 2 a_2 \geq 6 a_3$ and $n(n-2)a_n \geq (n-1)(n+1)a_{n+1}$ for $n \geq 3$. Then $f(z)=z+ \sum_{n=2}^\infty a_n z^n$ is close-to-convex with respect to both the starlike functions $z$ and $z/(1-z)$. Further, the function $f$ is starlike univalent in $\mathbb{D}$.
\end{lemma}

In our next result we will study the close-to-convexity of $z\mathcal{B}_\nu$ with respect to the starlike functions $z/(1-z^2)$. \begin{theorem}\label{thm:cc-2}
If $\nu \geq \nu_0 \approx 19.6203$, the function $z\mathcal{B}_\nu$ is  close-to-convex with respect to the starlike functions  $z/(1-z^2)$.
\end{theorem}
The following result is use to prove Theorem $\ref{thm:cc-2}$.
\begin{lemma}\cite[Theorem 4.4]{MS}\label{lemma-cc}
Let $\{a_n\}_{n\geq1}$ be a sequence of positive real numbers such that $a_1=1$. Suppose that $a_1 \geq 8 a_2$, and $(n-1)a_n \geq (n+1)a_{n+1}$ for all $n \geq 2$. Then the function $f(z)=z+ \sum_{n=2}^\infty a_n z^n$ is close-to-convex with respect to the starlike functions $z/(1-z^2)$.
\end{lemma}

\section{Proof of the main results}\label{sec-2}
\begin{proof}[\bf{\it Proof of Proposition \ref{prop1}}]
From \eqref{eqn:B-S-power}, it follows that
\begin{align}\label{eqn:B-S-bessel}
z^\nu \mathcal{B}_\nu(z)
&=\sum_{n=0}^\infty \frac{{\Gamma(\nu+1)}\Gamma{\left(\frac{n+1}{2}\right)}}{\sqrt{\pi} n! \Gamma\left(\frac{n}{2}+\nu+1\right)} z^{n+\nu}\\ \notag
&= \sum_{m=0}^\infty \frac{{\Gamma(\nu+1)}\Gamma{\left(m+\frac{1}{2}\right)}}{\sqrt{\pi} (2m)! \Gamma\left(m+\nu+1\right)} z^{2m+\nu}\\ \notag
&\quad \quad+\sum_{m=0}^\infty \frac{{\Gamma(\nu+1)}\Gamma{\left(m+1\right)}}{\sqrt{\pi} (2m+1)! \Gamma\left(m+\nu+\frac{3}{2}\right)} z^{2m+1+\nu}.
\end{align}
The Legendre  duplication formula (see \cite{ABRAMOWITZ,Andrews-Askey})
\[\mathrm{\Gamma}{(z)}\mathrm{\Gamma}{\left(z+\frac{1}{2}\right)}= 2^{1-2z}\; \sqrt{\pi}\; \mathrm{\Gamma}{(2z)}\]
shows that
\begin{align*}
\frac{\Gamma{\left(m+\frac{1}{2}\right)}}{\sqrt{\pi} (2m)!}= \frac{1}{2^{2m} m!} \quad \text{and}
\quad \frac{\Gamma{\left(m+1\right)}}{\sqrt{\pi} (2m+1)!}= \frac{1}{2^{2m+1} \Gamma{\left(m+\frac{3}{2}\right)}}.
\end{align*}
This along with \eqref{eqn:modified-Bessel} and \eqref{modStrv}, the identity \eqref{eqn:B-S-bessel} reduce to
\begin{align*}
z^\nu \mathcal{B}_\nu(z)
&= 2^\nu\Gamma(\nu+1)\sum_{m=0}^\infty \frac{\left(\frac{z}{2}\right)^{2m+\nu}}{m! \Gamma\left(m+\nu+1\right)} \left(\frac{z}{2}\right)^{2m+\nu}\\
&\quad +\sum_{m=0}^\infty \frac{\left(\frac{z}{2}\right)^{2m+\nu+1}}{\Gamma\left(m+\frac{3}{2}\right) \Gamma\left(m+\nu+\frac{3}{2}\right)}\\
&=2^\nu\Gamma(\nu+1) (I_\nu(z)+L_\nu(z)).
\end{align*}
This complete the proof.
\end{proof}

\begin{proof}[\bf{\it Proof of Proposition \ref{prop2}}]
Differentiating  the series  \eqref{eqn:B-S-power}, it follows that
\begin{align*}
z \frac{d}{dz}\mathcal{B}_\nu(z)
&= \sum_{n=0}^\infty \frac{n \Gamma{(\nu+1)} \Gamma{\left(\frac{n+1}{2}\right) }}{\sqrt{\pi} n! \; \Gamma{\left(\frac{n}{2}+\nu+1\right)}}z^n\\
&= 2\nu \sum_{n=0}^\infty \frac{\left(\frac{n}{2}+\nu\right) \Gamma{(\nu)} \Gamma{\left(\frac{n+1}{2}\right) }}{\sqrt{\pi} n! \; \Gamma{\left(\frac{n}{2}+\nu+1\right)}}z^n\\
&\quad -
2\nu\sum_{n=0}^\infty \frac{ \Gamma{(\nu+1)} \Gamma{\left(\frac{n+1}{2}\right) }}{\sqrt{\pi} n! \; \Gamma{\left(\frac{n}{2}+\nu+1\right)}}z^n\\
&= 2 \nu \mathcal{B}_{\nu-1}(z)-2\nu \mathcal{B}_\nu(z). \qedhere
\end{align*}
\end{proof}

\begin{proof}[\bf{\it Proof of Theorem \ref{thm:starlike}}]
Denote $f(z)= z\mathcal{B}_\nu(z)$. Then a computation together with the hypothesis \eqref{eqn:thm-starlike} yield
\begin{align*}
\left|\frac{zf'(z)}{f(z)}-1\right|
=\left|\frac{z\mathcal{B}'_\nu(z)}{\mathcal{B}_\nu(z)}\right|
<1-\lambda, \end{align*}
which is equivalent to \eqref{eqn:suf-starlike} for $\beta=0.$ The conclusion follows from Lemma \ref{lem:suf-starlike}.
\end{proof}

\begin{proof}[\bf{\it Proof of Theorem \ref{thm:starlike1}}]
Define $h : \mathbb{D} \to \mathcal{C}$ as
\begin{align}\label{eqn:h}
h(z):= \frac{2 \Gamma(\nu+3/2)}{\Gamma(\nu+1)} (\mathcal{B}_\nu(z)-1).
\end{align}
Then clearly $h \in \mathcal{A}$. Now a computation yield
\begin{align*}
\left|\frac{zh''(z)}{h'(z)}\right|= \left|\frac{z\mathcal{B}_\nu''(z)}{\mathcal{B}_\nu'(z)}\right| < \frac{2-3\nu+2\nu^2}{2(1-\nu)}.
\end{align*}
Taking $\beta=1$,  from Lemma $\ref{lem:suf-starlike}$, it follows that $h \in \mathcal{S}^{\ast}(\lambda)$ with respect to origin. Now Theorem \ref{thm:starlike1} follows from the definition of $h$  in \eqref{eqn:h}.
\end{proof}

\begin{proof}[\bf{\it Proof of Theorem \ref{theorem-1}}]
From  \eqref{eqn:B-S-power}, we can express $z \mathcal{B}_\nu$ as
\begin{align}
z \mathcal{B}_\nu(z) = \sum_{n=1}^\infty a_n z^n,
\end{align}
where
\begin{align}\label{eqn:coef}
a_1=1 \quad \text{and}  \quad a_n= \frac{ {\Gamma(\nu+1)}\Gamma{\left(\frac{n}{2}\right)}}{\sqrt{\pi} (n-1)!\; \Gamma\left(\frac{n+1}{2}+\nu\right)}, \quad  n \geq 2.
\end{align}
Define the function $g_n: [0, \infty) \to \mathbb{R}$ as
\begin{align}
g_n(\nu) := \frac{\Gamma\left(\frac{n}{2}+1+\nu\right)}{\Gamma\left(\frac{n+1}{2}+\nu\right)}.
\end{align}
The logarithmic differentiation with respect to $\nu$ implies
\begin{align}
g_n'(\nu)= g_n(\nu)\bigg( \Psi\left(\frac{n}{2}+\nu+1\right)- \Psi\left(\frac{n+1}{2}+\nu\right)\bigg).
\end{align}
Here $\Psi$ is the well-known digamma functions which is an increasing function on $[0, \infty)$,  and consequently $g_n$ is also increasing. Thus for $\nu \geq 1/2$,
\begin{align}
g_n(\nu) \geq g_n(1/2)=\frac{\Gamma\left(\frac{n+3}{2}\right)}{\Gamma\left(\frac{n}{2}+1\right)}
\end{align}
Thus for $n\geq 1$ and $\nu \geq 1/2$, it follows that
\begin{align}\label{eqn-ration-an}
\frac{a_n}{a_{n+1}}&= \frac{ {\Gamma(\nu+1)}\Gamma{\left(\frac{n}{2}\right)}}{\sqrt{\pi} (n-1)! \Gamma\left(\frac{n+1}{2}+\nu\right)}
\times \frac{\sqrt{\pi} (n)! \Gamma\left(\frac{n+2}{2}+\nu\right)} { {\Gamma(\nu
+1)}\Gamma{\left(\frac{n+1}{2}\right)}}\\ \notag
&\geq \frac{n \Gamma{\left(\frac{n}{2}\right)} \Gamma\left(\frac{n+3}{2}\right)}{\Gamma{\left(\frac{n+1}{2}\right)}\Gamma\left(\frac{n}{2}+1\right)}=n+1.
\end{align}
This implies  for $n\geq 1$
\begin{align*}
\underline{\Delta} a_n &= n a_n- (n+1) a_{n+1}\geq (n^2-1) a_{n+1} \geq 0,
\end{align*}
and for $n \geq 2$ we have
\begin{align*}
\underline{\Delta}^2 a_n&= n a_n- 2(n+1) a_{n+1}+(n+2) a_{n+2}\\
&\geq (n+1)(n-2)a_{n+1}+(n+2) a_{n+2}\\
&\geq (n+2) (n^2-n-1) a_{n+2}>0.
\end{align*}
 Thus $\{a_n\}$ satisfy the hypothesis of  Lemma $\ref{lemma-ac}$ and hence the conclusion.
\end{proof}

\begin{proof}[\bf{\it Proof of Theorem \ref{thm:cc-2}}]
The inequality \eqref{eqn-ration-an} yield that  for $n \geq 2$, and $\nu \geq 1/2$.
\begin{align*}
(n-1)a_n - (n+1)a_{n+1} \geq  n a_{n+1}>0,
\end{align*}
Now from \eqref{eqn:coef} it follows that the coefficient $a_n$ satisfy the hypothesis $a_1 \geq  8 a_2$ is equivalent to $\sqrt{\pi}\Gamma(\nu+3/2) \geq  8 \Gamma(\nu+1)$ which holds when $\nu \ge \nu_0$, where $\nu_0 \approx 19.6203$ is the positive root of the identity
\[ \sqrt{\pi}\Gamma(\nu+3/2) = 8 \Gamma(\nu+1).\]
Now the result follows from the Lemma $\ref{lemma-cc}$.
\end{proof}

\begin{problem} [Open]
Find the sharp lowest bound for $\nu >-1$ so that $z \mathcal{B}_\nu$ is starlike in $\mathbb{D}$ and also close-to-convex with respect to the starlike functions $z/(1-z^2)$.
\end{problem}

\section*{Acknowledgement} The author thanks the
Deanship of Scientific Research at King Faisal University for funding this work
under project number 150244.

\end{document}